\documentclass[a4paper,11pt]{article}
\usepackage{latexsym}
\usepackage{amssymb}
\usepackage{amsfonts}
\usepackage{amsmath}
\usepackage{indentfirst}
\usepackage{graphicx}
\usepackage{epsfig}
\newtheorem{prop}{Proposition}[section]
\newtheorem{defi}{Definition}[section]
\newtheorem{teo}{Theorem}[section]
\date{}

\author{Antonio Bernini \and Elisabetta Grazzini
\and Elisa Pergola \thanks{Communicating author: \emph{Phone:
+390554237458\ Fax: +390554237436\ e-mail:
elisa@dsi.unifi.it}.\newline \indent E-mails:
\texttt{\{bernini,ely,elisa,pinzani@dsi.unifi.it\}}} \and Renzo
Pinzani}

\title{A general exhaustive generation algorithm for Gray structures}

\begin{document}
\maketitle
\begin{center}
\vspace{-1cm} \noindent \scriptsize{Dipartimento di Sistemi e
Informatica, Universit\`a di Firenze. Viale G. B. Morgagni
65,\quad\quad\qquad 50134 Firenze, Italy.}
\end{center}

\vspace{1cm}

\begin{abstract}
Starting from a succession rule for Catalan numbers, we define a
procedure encoding and listing the objects enumerated by these
numbers such that two consecutive codes of the list differ only
for one digit. Gray code we obtain can be generalized to all the
succession rules with the \emph{stability property}: each label
$(k)$ has in its production two labels $c_1$ and $c_2$, always in
the same position, regardless of $k$. Because of this link, we
define \emph{Gray structures} the sets of those combinatorial
objects whose construction can be encoded by a succession rule
with the stability property. This property is a characteristic
that can be found among various succession rules, as the finite,
factorial or transcendental ones.

We also indicate an algorithm which is a very slight modification
of the Walsh's one, working in a $O(1)$ worst-case time per word
for generating Gray codes.
\end{abstract}

\section{Introduction} \label{intro}
The matter of encoding and listing the objects of a particular
class is common to several scientific topics, ranging from
computer science and hardware or software testing to chemistry,
biology and biochemistry.  Often, it is very useful to have a
procedure for listing or generating the objects in a particular
order. A very special kind of list is the so called Gray code,
where two successive objects are encoded in such a way that their
codes differ as little as possible (see below for more details and
\cite{Wa1}). There are many applications of the theory of Gray
codes for several combinatorial objects, involving permutations
\cite{J}, binary strings, Motzkin and Schr\"oder words \cite{V2},
derangements \cite{BV}, involutions \cite{Wa2}. They are also used
in other technological subjects as circuit testing, signal
encoding \cite{L}, data compression and other (we refer to
\cite{BBGP} for an exhaustive bibliography on the general matter).

The generation of a Gray code is often strictly connected with the
nature of the objects which we are dealing with. So, it seems to
have some importance the definition of a Gray code for the objects
of the classes with some common characteristic, as the classes
enumerated by the same sequence. From the idea of \cite{BBGP},
which we briefly recall in the sequel, in this work we develop a
procedure for listing the objects of those structures whose
exhaustive generation can be encoded by particular succession
rules (see below), say succession rules satisfying the
\emph{stability property} (see Section \ref{generalization}). In
order to point out the relation between such structures and the
possibility to list their objects in a Gray code, we define them
\emph{Gray structures}.

The starting point is the ECO method (see \cite{BDPP1} for a
survey). Closely related to this important enumerative tool is the
concept of \emph{succession rule} $\Omega$ \cite{CGHK,We1,We2},
which is a system defined by an \emph{axiom} $a \in \mathbb{N}$
and a set of \emph{productions}. The usual notation for a
succession rule is the following:
$$
\Omega: \left\{
\begin{array}{l}
   (a)\\
   (k)\rightsquigarrow(e_1(k))(e_2(k))\ldots(e_k(k))\ ,\quad \quad k\in \mathbb{N}\\
\end{array}
\right.
$$
The succession rule $\Omega$ can also be described with a rooted
tree where the nodes are the labels of $\Omega$: the axiom $(a)$
is the root of the tree and each node with label $(k)$ generates
$k$ sons with labels $(e_1(k)),(e_2(k)),\ldots,(e_k(k))$. The
structure we obtain is the so called  \emph{generating tree} of
$\Omega$ \cite{BETAL,CGHK}.

Our discussion moves from the well known succession rule
$\Omega_C$,
$$\Omega_C: \left\{
\begin{array}{l}
  (2) \\
  (k)\rightsquigarrow(2)(3)\ldots(k)(k+1),\quad \quad k\geq 2 \\
\end{array}
\right.
$$
defining the sequence of Catalan numbers and whose first levels of
the related generating tree are shown in Figure \ref{fig:albero e
generazione}.
Each object $x$ with size $n$ corresponds to a node at level $n-1$
(being the root of the tree at level $0$, corresponding to the
object of size $1$) and can be described by a word $w=w_1w_2\ldots
w_n$ encoding the path from the root to the node corresponding to
$x$: each $w_i$ is the label of a node of the path and is
generated by $\Omega_C$. In \cite{BBGP} the authors give a method
to exhaustively generate all the objects (words) of a given size
$n$ which substantially coincide with the reading from left to
right in the $(n-1)$-th level of the tree. So, the words at level
3 are generated in the following order (see figure \ref{fig:albero
e generazione}): \small
$$2222, 2223,
2232, 2233, 2234, 2322, 2323, 2332, 2333, 2334,2342, 2343, 2344,
2345.$$ \normalsize In the above list it is possible that two
consecutive words differ more than one digits: for instance, 2223
and 2232 differing in two digits or 2234 and 2322 with three
different digits. Our aim is to generate all the words of length
$n$ (naturally without repetitions) in such a way that \emph{two
consecutive words differ only for one digit}. Such a property is
strictly related to the concept of ``Gray code", which definition
we relate can be found in \cite{Wa1}. Substantially, it can be
summed up in the following: \emph{a Gray code is an infinite set
of word-lists with unbounded word-length such that the Hamming
distance between any two adjacent words is bounded independently
of the word-length} (the Hamming distance is the number of
positions in which two words differ). For a complete discussion on
Gray codes we refer the reader to the paper of T. Walsh
\cite{Wa1}.

In Section \ref{idea} an informal description of the used strategy
for our purpose is presented, referring to objects whose
construction can be described by $\Omega_C$. Then, in Section
\ref{definition}, a rigorous definition of the list (Definition
\ref{lista}) and a proof that it is a Gray code (Theorem
\ref{teoGray}) are given. Section \ref{Dyck} presents the
application and the analysis to the particular case of Dyck paths,
enumerated by Catalan numbers. Finally, Section
\ref{generalization} generalizes the construction of the Gray code
to those objects whose generation can be described by succession
rules with the \emph{stability property}. In that section, we also
present some examples of Gray structures.

\section{The procedure}\label{idea}
The strategy used in \cite{BBGP} for listing the objects of size
$n$ corresponds to a visit of all the nodes at level $n-1$ in the
generating tree from left to right. So, after the visit of a
subtree $T_i$ is completed, the path from the root to the leftmost
node of the successive subtree $T_{i+1}$ has at least two
different nodes with respect to those ones of the last path of the
preceding subtree $T_i$. This is due to the fact that the labels
of the sons of a node are visited in the same order they have in
the production of the succession rule $\Omega_C$, where the list
of the successors of a label $(k)$ is $<2, 3, \ldots, k, k+1>$.

For our purpose we must check that when a subtree has been
completely visited and if $v$ is the last path generated in such a
visit, then the successive path $w$ has only one different digit
with respect to the digits of $v$. We now illustrate the procedure
we are going to use referring to Figure \ref{fig:albero e
generazione}, where the words of length 4 are generated.
\begin{figure}
\begin{center}
\includegraphics[scale=.5]{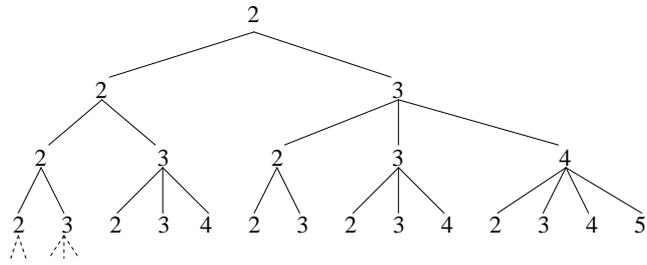}
\includegraphics[scale=0.5]{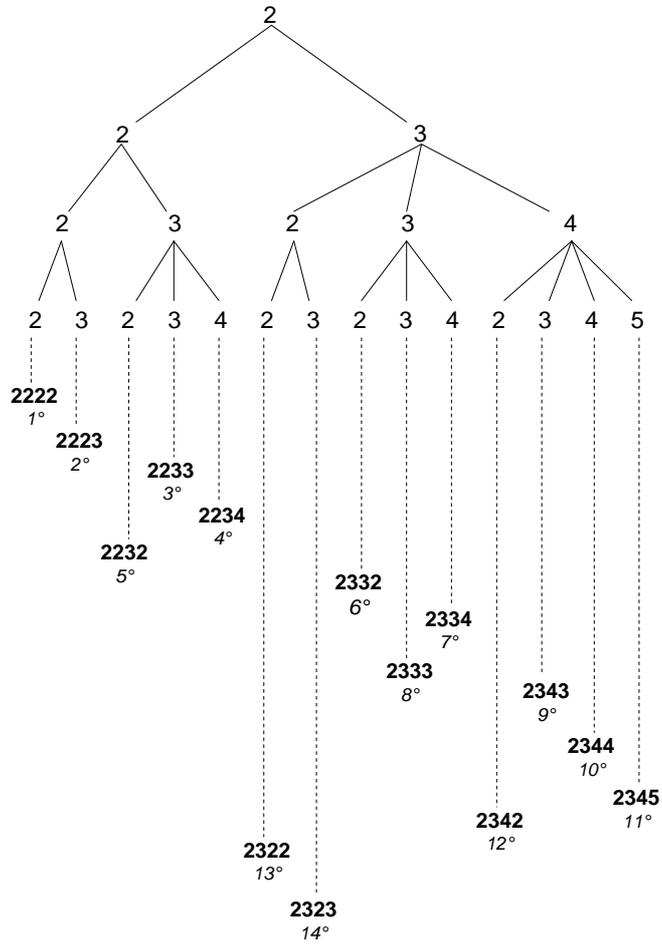}
\end{center}
\caption{first levels of the generating tree for Catalan numbers
(upper figure); generation of the words of length 4 (lower
figure).}\label{fig:albero e generazione}
\end{figure}

The first object of the list is the word 2222, corresponding to
the path from the root to the leftmost node at level 3 in the
generating tree. Then, in order to complete the visit of the
current subtree, the second word is 2223. At this point, the next
path in the list will have a different digit with respect to the
digits of 2223, which is not the last one: in order to respect the
above definition of Gray code, the third word in the list could be
2323 or 2233. The choice is determined by the leading idea that a
successive path $w=aw_2\ldots w_n$ must have as much as possible
the same edges of the preceding path $v=av_2\ldots v_n$ in the
list and if $v_j$ and $w_j$ are the first nodes necessarily
different in $v$ and $w$, then all the nodes $v_r$ and $w_r$ must
have the same labels for $r=j+1,\ldots,n-1,n$, in order to respect
the Gray code definition. So, the third word is 2233. The fourth
and the fifth one are 2234 and 2232, respectively. From the
generation of these last two words we can deduce that only the
last digit is changed when a same subtree is visited and that the
order for changing the last digit is \emph{shifted} with respect
to the classical one in a cyclic way in order to complete the set
of the sons of the second-last digit: for the sake of clearness in
this case the shifted list of the successors of the second-last
digit 3 is $<3, 4, 2>$, while the classical one would be $<2, 3,
4>$. This fact can be generalized. Let $e$ be the first path of a
new subtree and let $i$ and $k$ be the the last and the
second-last digit of $e$, respectively ($i\neq 2$, see below).
Then the right order for changing the last digit is
$<i,i+1,\ldots,k,k+1,2,3,\ldots,i-1>$.

The sixth path which is now generated is $f=2332$, according to
the above leading idea. Note that the second digit is changed with
respect of the second digit of the fifth word and that the third
and the fourth digits in $f$ are the same you find in 2232. The
word $f$ is the first path of a new subtree and then only the last
digit has to be changed, till the whole set of the sons of the
second-last digit 3 is completed. Since the last digit of $f$ is
2, one could think that in this case the \emph{shifted production}
of the digit 3 coincides with the classical production $<2,3,4>$,
obtaining that the 6-th, 7-th and 8-th words are 2332, 2333, 2334,
respectively. But so doing the procedure fails when it is used to
list the words of length 6, as the reader can easily check when he
arrives at the generation of the word 234565. The reason of the
failure will be clear in the next section, where the rigorous
formalization of our procedure is presented. The right way for
changing the last digit of $f$ is to consider the list $<2,4,3>$
of the sons of the digit 3, then obtaining the 6-th, 7-th and 8-th
words as follows: 2332, 2334 and 2333, respectively. This fact
suggest us that if $f$ is the first word of a new subtree, if its
last digit is 2 and if $k$ is its second-last digit, then the
right order for changing the last digit is
$<2,k+1,k,k-1,k-2,\ldots,4,3>$. The remaining objects can be now
easily obtained, as in Figure \ref{fig:albero e generazione}.

We now summarize the definition of the \emph{shifted production}
which is used to change the last digit in the words. Let
$v=v_1v_2\ldots v_n$ be the first path of a new subtree. Let $k$
and $i$ be the second-last and the last digit of $v$,
respectively, then the list $s(k,i)$ of the sons of $k$ such that
the first son is $i$, is:
$$\left\{
\begin{array}{l}
  s(k,2)=<2,k+1,k,k-1,\ldots 4,3> \\
  s(k,i)=<i,i+1,\ldots, k-1,k,k+1,2,3,\ldots,i-1> \quad.\\
\end{array}
\right.
$$

\section{A Gray code for Catalan structures}\label{definition}
First we define the lists for the objects whose generating tree
can be described by the succession rule for the Catalan numbers we
presented in the previous section, then we will prove (Theorem
\ref{teoGray}) that these lists form a Gray code, in the sense of
the definition in Section \ref{intro}. The following notation is
used:
\begin{itemize}
    \item $\mathcal{L}_k=$ list of the codes of the objects of
    length $k$;
    \item $l^k_i=$ $i$-th element of $\mathcal{L}_k$;
    \item $|\mathcal{L}_k|$= cardinality of $\mathcal{L}_k$;
    \item if $x$ is a sequence of digits, then $\overrightarrow{x}$ is the
    rightmost digit of $x$;
    \item $\Theta$ is the concatenation of lists;
    \item if $L$ is a list, then:
        \begin{itemize}
            \item \emph{first}($L$) denotes the first element of the
            list $L$;
            \item \emph{last}($L$) denotes the last element of
            the list $L$;
            \item $x\circ L$ is the list obtained by pasting $x$ with
            each element of $L$.
        \end{itemize}
\end{itemize}
Our definition is a recursive definition and it is based on a
generation of sublists with increasing length:
\begin{defi}\label{lista}
The list $\mathcal{L}_n$ of all the elements of length $n$ is
$$
\left\{
\begin{array}{cccc}
  \mathcal{L}_1 & = & <2>& \\
  & & &\\
  \mathcal{L}_n & = & \Theta_{i=1}^M\ L_n^i&\quad \quad \mbox{if}\quad n>1\\
\end{array}%
 \right.
 $$
 where $M=|\mathcal{L}_{n-1}|$ and $L_n^i$ is defined by
 $$
 \left\{
 \begin{array}{cccc}
    L_n^1&=&l_1^{n-1}\circ s(2,2)&\\
    & & &\\
    L_n^i&=&l_i^{n-1}\circ s(\overrightarrow{l_i^{n-1}},\overrightarrow{last(L_n^{i-1})})&\mbox{if}\quad i>1 \quad \quad .\\
 \end{array}
  \right.
 $$
\end{defi}
The list $L_n^1$ is obtained by linking together the first element
of the list of the objects of size $n-1$ (i.e. $l_1^{n-1}$) with
the elements of the list $s(2,2)=\ <2,3>$; then $L_n^1$ has always
two elements: $l_1^{n-1}2$ and $l_1^{n-1}3$. The next lists
$L_n^i$ with $i>1$ are obtained as follows:
\begin{itemize}
    \item consider the $i$-th element of $\mathcal{L}_{n-1}$ (i.e.
    $l_i^{n-1}$);
    \item consider the list of the successors of the rightmost
    digit of $l_i^{n-1}$ shifted starting from the rightmost
    digit of the rightmost element of $L_n^{i-1}$ (i.e.
    $s(\overrightarrow{l_i^{n-1}},\overrightarrow{last(L_n^{i-1})})$);
    \item paste $l_i^{n-1}$ with each element of the list
    $s(\overrightarrow{l_i^{n-1}},\overrightarrow{last(L_n^{i-1})})$.
\end{itemize}
Let us construct for instance the list $\mathcal{L}_4$:
\begin{description}
    \item \begin{description}
                \item $\mathcal{L}_1=<2>;$
          \end{description}
    \item $L_2^1=2 \circ s(2,2)=2 \circ <2,3>=<22,23>$, then
    \item \begin{description}
                \item $\mathcal{L}_2=<22,23>;$
          \end{description}
    \item $L_3^1=22\circ s(2,2)=22\circ <2,3>=<222,223>$;
    \item $L_3^2=23\circ s(3,3)=23\circ <3,4,2>=<233,234,232>$, then
    \item \begin{description}
                \item $\mathcal{L}_3=<222,223,233,234,232>;$
          \end{description}
    \item $L_4^1=222\circ s(2,2)=222\circ <2,3>=<2222,2223>$;
    \item $L_4^2=223\circ s(3,3)=223\circ <3,4,2>=<2233,2234,2232>$;
    \item $L_4^3=233\circ s(3,2)=233\circ <2,4,3>=<2332,2334,2333>$;
    \item $L_4^4=234\circ s(4,3)=234\circ <3,4,5,2>=<2343,2344,2345,2342>$;
    \item $L_4^5=232\circ s(2,2)=232\circ <2,3>=<2322,2323>$, then
    \item \begin{description}
                \item
                \begin{eqnarray*}
                \mathcal{L}_4 &=& <2222,2223,2233,2234,2232,2332,2334,2333,2343,\\
                & &2344,2345,2342,2322,2323>
                \end{eqnarray*}
          \end{description}
\end{description} We now prove the following:
\begin{teo} \label{teoGray}
Two consecutive elements of the list $\mathcal{L}_n$ differ only
for one digit.
\end{teo}

\emph{Proof.} We can proceed by induction on $n$:
\begin{description}
    \item[\textbf{base:}] if $n=1$, then the theorem is trivially
    true since $\mathcal{L}_1=<2>$;
    \item[\textbf{inductive hypothesis:}] let us suppose that
    $l_i^{n-1}$ and $l_{i+1}^{n-1}$, with $1\leq i\leq |\mathcal{L}_{n-1}|-1$, differ only for one
    digit;
    \item[\textbf{inductive step:}] the list $\mathcal{L}_n$ is
    obtained by linking together the lists $L_n^i$ for $i=1,\ldots
    ,|\mathcal{L}_{n-1}|$. Since the elements of each list $L_n^i$
    differ only for one digit by construction, we must prove
    the statement only for $last(L_n^i)$ and $first(L_n^{i+1})$,
    with $1\leq i\leq |\mathcal{L}_{n-1}|-1$.

    Let $J$ be the last element of
    $s(\overrightarrow{l_i^{n-1}},\overrightarrow{last(L_n^{i-1})})$.
    Then we have:
    $$
    last(L_n^i)=l_i^{n-1}J\ .
    $$
    We also have:
    $$
    L_n^{i+1}=l_{i+1}^{n-1}\circ
    s(\overrightarrow{l_{i+1}^{n-1}},\overrightarrow{last(L_n^i)})=l_{i+1}^{n-1}\circ
    s(\overrightarrow{l_{i+1}^{n-1}},J)\ .
    $$
    From the definition of the shifted list of the successors we
    deduce that the first element of a list $s(i,k)$ is always
    $k$, then:
    $$
    first(L_n^{i+1})=l_{i+1}^{n-1}J\ .
    $$
    Since $l_i^{n-1}$ and $l_{i+1}^{n-1}$ differ only for one
    digit by the inductive step, this statement holds also for
    $last(L_n^i)$ and $first(L_n^{i+1})$. So, the theorem is
    proved. \begin{flushright}$\square$\end{flushright}\phantom\\

\end{description}

At this point it is easily seen that
$\overrightarrow{first(L_n^{i+1})}$ is a son of the second-last
digit of $first(L_n^{i+1})$ and that
$\overrightarrow{first(L_n^{i+1})}=\overrightarrow{last(L_n^i)}$.
We remark that it is not possible that
$\overrightarrow{last(L_n^i)}$ does not belong to the set of sons
of the second-last digit of $first(L_n^{i+1})$, since from the
definition of the \emph{shifted production}, the construction we
described above and the axiom of $\Omega_C$ (which is 2), we
deduce that $\overrightarrow{last(L_n^i)}\in \{2,3\}$, which are
present in the production of each possible label.


\subsection{The algorithm to generate $\mathcal{L}_n$}

The aim is defining an algorithm which is not recursive for
generating all the words of length $n$ encoding the objects of
size $n$. We base our procedure on the general idea that if a word
$c_j$ has been generated, then a single digit must be changed to
generate the next word $c_{j+1}$, as the authors made in
\cite{BBGP}.

The first word of the list is $w=222\ldots 2$, where $w_i=2$, for
$i=1,2,\ldots n$. The digit $w_i$ to be modified at each step is
determined using the algorithm of Walsh \cite{Wa1}, i.e. using a
$(n+1)$-dimensional array $e$, which is updated in such a way
that, at each step, $e_{n+1}$ points to $w_i$. Once $w_i$ is
determined, it can not be modified by simply increasing it by one
\cite{BBGP}, but the definition of the \emph{shifted production}
must be taken in account. So, we use another array $d$
($n$-dimensional), which is defined as follows: $d_i=0$ if $w_i$
is modified according to the shifted production $s(k,2)$; $d_i=1$
if $w_i$ is modified according to $s(k,3)$. It is easy to prove
that the introduction of the array $d$ does not exchange the
complexity of the recalled procedure of Walsh for generating Gray
codes in $O(1)$ worst-case time per word \cite{Wa1}: his clever
algorithm remains the starting point for the implementation of our
method.

We note that $d$ can also be used to establish when $w_i$ is no
more modifiable: from the definition of $s(k,j)$ it happens if
$(d_i=0 \wedge w_i=3)$ or if $(d_i=1 \wedge w_i=2).$

The generating procedure stops when the digit to be modified is
$w_1$.

\section{The case of Dyck paths}\label{Dyck}
We consider now the specific class of Dyck paths. Each of them can
be associated with a binary string according to the substitution,
for example, of the up steps with the $1$ bit and the down steps
with $0$. Let us consider a word of length $n$ of the Gray code
defined in Section \ref{definition}. It has a correspondent Dyck
path which, in turn, is associated with a binary string, both of
length $2n$ (in Section \ref{gray-bin} we present an algorithm to
directly translate a word in the associated binary string). We
want to prove that, if we consider two consecutive binary strings
corresponding to two consecutive words in the Gray code, they
differ only for two bits (note that the \emph{Hamming distance}
between two binary strings encoding two Dyck paths is at least
$2$). For this aim we base on the ECO construction of Dyck paths
\cite{BDPP1}. We recall briefly its main features: if $p$ is a
Dyck paths of length $2n$ with the last descent of $k$ steps, then
it has $k+1$ active sites; we obtain each of its $k+1$ sons by
inserting a peak in each active sites; the insertion of a peak in
an active sites at hight $h$ generates a Dyck path with $h+2$
active sites.
Now we state the next proposition:
\begin{prop}\label{teoremalungo}Two words of the Gray
code differing for one digit correspond to binary strings which
differ only for two bits.
\end{prop}

\emph{Proof.} The last digit of a word denotes the number of
active sites of the corresponding Dyck path, so if it is $k$, then
the path has $k-1$ down steps in the last descent, according to
the above mentioned ECO construction.
\begin{description}
    \item[A] Let us consider the case when the two words differs in
    the last digit. Let their codes be:
    $$w_1w_2\ldots w_iw_{i+1}$$ and $$w_1w_2\ldots w_iz_{i+1}.$$
    We indicate a generic bit with the star $*$, so $w_1\ldots
    w_i$ corresponds to
    $$
    \underbrace{1**\ldots**1}_{2i-w_i+1}\underbrace{000\ldots
    0}_{w_i-1}.
    $$
    The adding of $w_{i+1}$ corresponds to the insertion of a peak
    at height $w_{i+1}-2$ in the last descent of the Dyck path associated to $w_1w_2\ldots w_i$.
    So, the corresponding binary string is
        \begin{equation}\label{pippo}
    \small{\underbrace{1**\ldots**1}_{2i-w_i+1}\underbrace{00\ldots\ldots\ldots 0}_{(w_i-1)-(w_{i+1}-2)}
    1\underbrace{00\ldots0}_{w_{i+1}-1}
    =
    \underbrace{1**\ldots**1}_{2i-w_i+1}\underbrace{00\ldots0}_{w_i-w_{i+1}+1}
    1\underbrace{00\ldots0}_{w_{i+1}-1}}
        \end{equation}
    (note that after the adding of $w_{i+1}$, the total number of bits is properly
    $2i+2$). In particular we have:
        \begin{itemize}
    \item in the case $w_{i+1}=w_i+1$, when the peak is inserted in the active site with maximal height,
    the binary string becomes
    $$
    \underbrace{1**\ldots **1}_{2i-w_i+1}1\underbrace{00\ldots
    0}_{w_i},
    $$
    in other words, the last ascent is longer than one step with respect to the Dyck path codified
    by the word $w_1w_2\ldots w_i$;
    \item in the case $w_{i+1}=2$, when the peak is added at
    height 0 at the end of the Dyck path corresponding to $w_1w_2\dots
    w_i$,the binary string is
    $$
    \underbrace{1**\ldots **1}_{2i-w_i+1}\underbrace{00\ldots
    0}_{w_i-1}10.
    $$
        \end{itemize}
    In a similar manner, the addition of $z_{i+1}$ after $w_i$
    transforms the corresponding binary string in
    $$
    \underbrace{1**\ldots**1}_{2i-w_i+1}\underbrace{00\ldots0}_{w_i-z_{i+1}+1}
    1\underbrace{00\ldots0}_{z_{i+1}-1}.
    $$
    Let us suppose that $z_{i+1}=w_{i+1}+j$, where $j$ can also
    assume negative values. If $j>0$, then $j\in \{1, w_i+1-w_{i+1}\}$;
    if $j<0$, then $j\in \{-1, 2-w_{i+1}\}$. The word $w_1w_2\ldots
    w_iz_{i+1}$ corresponds to the binary string
        \begin{equation}\label{pluto}
    \underbrace{1**\ldots**1}_{2i-w_i+1}\underbrace{00\ldots0}_{w_i-w_{i+1}+1-j}
    1\underbrace{00\ldots0}_{w_{i+1}-1+j}
        \end{equation}
The difference between the words (\ref{pippo}) and (\ref{pluto})
is the location of the rightmost $1$ bit, which in (\ref{pluto})
is shifted of $|j|$ positions towards left $(j>0)$ or right
$(j<0)$ with respect to (\ref{pippo}). It easily seen that the two
strings differ only for the two bits in position $w_{i+1}$ and
$w_{i+1}+j$ from the right of the word.
    \item[B] Let us consider now the case when the two words differ
    for two digits which are not the last ones:
    \begin{equation}\label{par1}
    w_1w_2\ldots w_iw_{i+1}w_{i+2}\ldots w_n
    \end{equation}
    and
    \begin{equation}\label{par2}
    w_1w_2\ldots w_iz_{i+1}w_{i+2}\ldots w_n.
    \end{equation}
    The associated binary strings after the insertion of $w_{i+2}$
    (i.e. the binary strings coding $w_1\ldots w_iw_{i+1}w_{i+2}$
    and $w_1\ldots w_iz_{i+1}w_{i+2}$) are
    $$
    \underbrace{1**\ldots**1}_{2i-w_i+1}\underbrace{00\ldots0}_{w_i-w_{i+1}+1}
    1\underbrace{00\ldots0}_{w_{i+1}-w_{i+2}+1}1\underbrace{00\ldots0}_{w_{i+2}-1}
    $$
    and
    $$
    \underbrace{1**\ldots**1}_{2i-w_i+1}\underbrace{00\ldots0}_{w_i-w_{i+1}+1-j}
    1\underbrace{00\ldots0}_{w_{i+1}-w_{i+2}+1+j}1\underbrace{00\ldots0}_{w_{i+2}-1}.
    $$
    where, as in the preceding case, $z_{i+1}=w_{i+1}+j$.
    The insertions of the next digits $w_k$ with $k=i+3,\ldots,n$, which are
    equal in the two words, modify in the same way the last descent in the associated Dyck
    paths. Then, the difference between the two binary strings
    corresponding to them is not due to these insertions. So, also in this case, the binary
    strings related to (\ref{par1}) and (\ref{par2}) differ only for
    two bits. \begin{flushright}$\square$\end{flushright}\phantom\\
\end{description}

\subsection{From a binary string to the next one}\label{next one
string}

The structure of the above proof can be used to derive an
algorithm to generate a binary string $p_{h+1}$ from the preceding
one $p_h$, taking into account the generation order of the
corresponding words in the Gray code. If $u_h$ and $u_{h+1}$ are
two consecutive words in the Gray codes and $p_h$ is the binary
string corresponding to $u_h$, then:
\begin{itemize}
    \item if $u_h$ and $u_{h+1}$ differ in the last digit and
    $j=\overrightarrow{u_{h+1}}-\overrightarrow{u_h}$
    is the difference between these ones, then $p_{h+1}$ is
    obtained from $p_h$ by the shifting of $|j|$ positions of the rightmost 1 bit towards
    left if $j>0$ or right if $j<0$;
    \item if $u_h$ and $u_{h+1}$ differ in the $i$-th digit
    and $j$ is the difference between the $i$-th digit of $u_{h+1}$ and the $i$-th
    digit of $u_h$, then $p_{h+1}$ is obtained from
    $p_h$ by the shifting of $|j|$ positions of the second rightmost 1 bit
    towards left if $j>0$ or right if $j<0$.
\end{itemize}
The correctness of the above procedure can be easily checked and
the algorithm is based on the proof of the preceding proposition.

\subsection{From the word to the binary string}\label{gray-bin}
The proof of Proposition \ref{teoremalungo} suggests also the idea
for an inductive algorithm which allows to derive the binary
string corresponding to a given word in the Gray code. Let us
suppose we have already encoded a word $w_1\ldots w_{n-1}$ in the
binary string $u$. The adding of a new digit $w_n$ modifies only
the final part of $u$, as we can deduce from the first part of the
proof of Proposition \ref{teoremalungo}. In particular, the
$w_{n-1}-1$ rightmost $0$ bits of $u$ corresponding to the last
descent of the related Dyck path, are replaced by $w_{n-1}+1$ bits
as in the following:
$$
\underbrace{000\ldots 0}_{w_{n-1}-1}\longrightarrow \underbrace{
\underbrace{000\ldots 0}_{w_{n-1}-w_n+1}1\underbrace{000\ldots
0}_{w_n-1}}_{w_{n-1}+1}
$$
It correspond to the adding of a peak in some site of the last
descent of the Dyck path related to $u$.

Then, starting from the binary string $10$ encoding the minimal
Dyck path whose relating word in the Gray code is $2$, it is
possible to get the binary string corresponding to $w_1\ldots w_n$
from the knowledge of that one related to $w_1\ldots w_{n-1}$ by
means of the following inductive procedure:

\begin{description}
    \item[\textbf{base:}]the binary string corresponding to the
    word $2$ is $10$;
    \item[\textbf{inductive hypothesis:}]assume that $u$ is the binary string
    codifying $w_1\ldots w_{n-1}$;
    \item[\textbf{inductive step:}]then the binary string
    corresponding to ${w_1\ldots w_n}$ is obtained replacing the $w_{n-1}-1$
    rightmost $0$ bits of $u$ with the $w_{n-1}+1$ bits $\underbrace{000\ldots 0}_{w_{n-1}-w_n+1}1
    \underbrace{000\ldots 0}_{w_n-1}$\ \ .
\end{description}
In the following example  the encoding of the word 2334 is shown:
$$
\begin{array}{ccccccc}
  10 & \rightarrow & 1100 & \rightarrow & 110100 & \rightarrow & 11011000 \\
  (2) &  & (23) &  & (233) &  & (2334)\\
\end{array}
$$
 \begin{flushright}$\square$\end{flushright}\phantom\\

\textbf{Note.} The algorithm of Section \ref{next one string}
allows to find a binary string $p_{h+1}$ starting from the
preceding one $p_j$ and the words $u_h$ and $u_{h+1}$ of the Gray
code, corresponding to $p_h$ and $p_{h+1}$, respectively. The
algorithm of this section, whereas, generates the binary string
from the corresponding word by means of an inductive procedure
which can turn out too heavy for large values of $n$ (the length
of the word).

Hence, the preceding algorithm, having a low complexity, can be
used to generate $p_{h+1}$ in the case the string $p_h$ and the
words $u_h$ and $u_{h+1}$ are known.

\section{Generalization to stable succession rules}\label{generalization}
The crucial point in the construction of the lists $\mathcal{L}_n$
is that each label $k$ in the succession rule $\Omega_C$ has in
its production the two labels $2$ and $3$, as we pointed out at
the end of Section \ref{definition}. This property, together with
the definition of the \emph{shifted production} of $k$, allows
$last(L_n^i)$ and $first(L_n^{i+1})$ to be different only for one
digit (which is not the last one). Starting from this remark, we
generalize the procedure to define the Gray code to all those
succession rules having a particularity similar to $\Omega_C$
which we would like to call \emph{stability property}, meaning
with this name that in each production of $k$ we always find two
labels, say $c_1$ and $c_2$, regardless of $k$.

\begin{defi}(stability property) We say that a succession rule $\Omega$
$$
\Omega: \left\{
\begin{array}{l}
   (a)\\
   (k)\rightsquigarrow(e_1(k))(e_2(k))\ldots(e_k(k))\ ,\quad \quad k\in \mathbb{N}\\
\end{array}
\right.
$$
is \emph{stable} if for each $k$ there exist two indexes $i,j$
($i<j$) such that $e_i(k)=c_1$ and $e_j(k)=c_2$ ($c_1\leq c_2$).
\end{defi}

\begin{figure}
\begin{center}
\includegraphics[scale=.7]{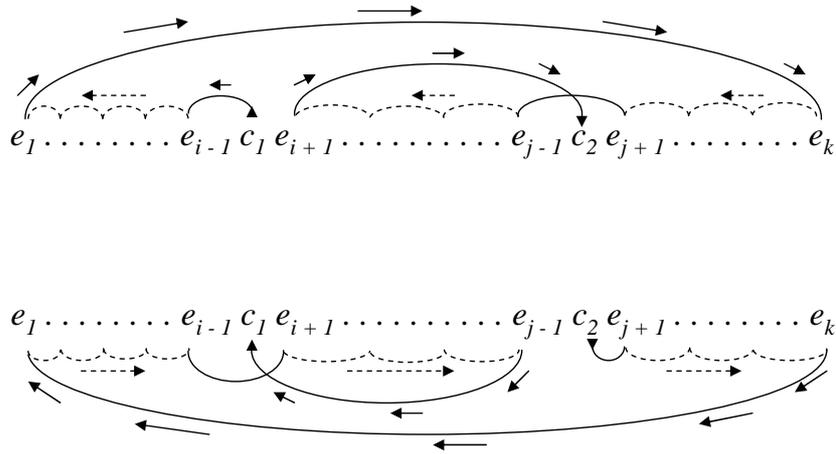}
\end{center}
\caption{Generalized shifted production.}\label{fig:figproduz}
\end{figure}

We need also to extend the definition of shifted production for
the labels of succession rules with the stability property, in
order to obtain that each list of successors of any $k$ ends with
$c_1$ or $c_2$. We have the following \emph{generalized shifted
productions} of $k$, being $e_i(k)=c_1$ and $e_j(k)=c_2$: \small
$$
\left\{
\begin{array}{l}
   s(k,c_1)=<c_1,e_{i-1}(k),\ldots,e_1(k),e_k(k),\dots,e_{j+1}(k),e_{j-1}(k),\ldots,e_{i+1}(k),c_2>\\
   \\
   s(k,c_2)=<c_2,e_{j+1}(k),\ldots,e_k(k),e_1(k),\ldots,e_{i-1}(k),e_{i+1}(k),\ldots,e_{j-1}(k),c_1>\ .\\
\end{array}
\right.
$$
\normalsize
In Figure \ref{fig:figproduz} we used two walks, very
similar to the \emph{factorial walks} on the integer half-line
\cite{BETAL}, to visualize the generalized shifted production of
$k$, the above one starting from $c_1$ and ending in $c_2$
(corresponding to $s(k,c_1)$) and  the below one starting from
$c_2$ and ending in $c_1$ (corresponding to $s(k,c_2)$).

Now, it is easy to prove that:

\begin{prop} If $\Omega$ is a succession rule with the stability
property, then the lists $\mathcal{L}_n$ defined by:
$$
\left\{
\begin{array}{cccc}
  \mathcal{L}_1 & = & <a>& \\
  & & &\\
  \mathcal{L}_n & = & \Theta_{i=1}^M\ L_n^i&\quad \quad \mbox{if}\quad n>1\\
\end{array}%
 \right.
 $$
 where $M=|\mathcal{L}_{n-1}|$ and $L_n^i$ is defined by
$$
  \left\{
 \begin{array}{cccc}
    L_n^1&=&l_1^{n-1}\circ s(\overrightarrow{l_1^{n-1}},c_1)&\\
    & & &\\
    L_n^i&=&l_i^{n-1}\circ s(\overrightarrow{l_i^{n-1}},\overrightarrow{last(L_n^{i-1})})&\mbox{if}\quad i>1\\
 \end{array}
  \right.
 $$
form a Gray code in the sense of the definition of Section
\ref{intro}, where two consecutive words of length $n$ differ for
one digit (Hamming distance equals to one).
\end{prop}

The proof is completely similar to that one of Theorem
\ref{teoGray} and it is omitted.


\vspace{.3cm}

Note that in the special case $i=1,\ j=2$ the generalized shifted
production is:
$$
\left\{
\begin{array}{l}
   s(k,c_1)=<c_1,e_k(k),e_{k-1}(k),\ldots,e_3(k),c_2>\\
   \\
   s(k,c_2)=<c_2,e_3(k),\ldots,e_k(k),c_1>\ .\\
\end{array}
\right.
$$

\vspace{1cm} We now analyze some particular cases of succession
rules with the stability property.

\vspace{.3cm} \noindent \textbf{Example 1.} Let us consider the
following rule $\Omega_{F_o}$,
$$
\Omega_{F_o}: \left\{
\begin{array}{l}
    (2)\\
    (2)\rightsquigarrow(2)(3)\\
    (3)\rightsquigarrow(2)(3)(3)\ ,
\end{array}
\right.
$$
defining the odd Fibonacci numbers. It is easily seen that it
satisfies the stability property, but the rule $\Omega_F$,
$$
\Omega_F: \left\{
\begin{array}{l}
    (2)\\
    (2)\rightsquigarrow(1)(2)\\
    (1)\rightsquigarrow(2)\ ,
\end{array}
\right.
$$
defining Fibonacci numbers, does not satisfy the stability
property. This is to say that such a property is not common to all
the succession rules of a certain family (finite succession rules,
in this case).

In the following examples it is shown that a similar behavior can
be found also in factorial or transcendental
rules. \begin{flushright}$\square$\end{flushright}\phantom\\

\vspace{.3cm}\noindent \textbf{Example 2.} The factorial rule:
$$
\Omega_M: \left\{
\begin{array}{l}
  (1) \\
  (k)\rightsquigarrow(1)(2)\ldots(k-1)(k+1)\ \ \ , \\
\end{array}
\right.
$$
defining the sequence of Motzkin numbers, does not satisfies the
stability property, since only for $k\geq 3$ each label has
$c_1=1$ and $c_2=2$ in its production. But the rules $\Omega_A$ of
kind
$$
\Omega_A: \left\{
\begin{array}{l}
   (a)\\
   (k)\rightsquigarrow(a)(a+1)\ldots (k)(k+1)(k+d_1)\ldots(k+d_m)\\
\end{array}
\right.
$$
(with $a\geq 2$, $m=a-2$, $d_i\geq 0$ and $d_i\leq d_{i+1}$) are
factorial and stable rules, with $i=1$, $j=2$, $c_1=a$ and
$c_2=a+1$. The following well-known succession rule $\Omega_t$,
related to the Gray structure of the $t$-ary trees \cite{BDP}, is
a particular case:
$$
\Omega_t: \left\{
\begin{array}{l}
   (t)\\
   (k)\rightsquigarrow(t)(t+1)\ldots(k-1)(k)(k+1)\ldots(k+t-2)(k+t-1)\\
\end{array}
\right.
$$
and the generalized shifted production is:
$$
\left\{
\begin{array}{l}
   s(k,t)=<t,k+t-1,k+t-2,\ldots,k+1,k,k-1,\ldots,t+2,t+1>\\
   \\
   s(k,t+1)=<t+1,t+2,\ldots,k-1,k,k+1,\ldots,k+t-2,k+t-1,t>\ .\\
\end{array}
\right.
$$
In the following, we present the construction of the list
$\mathcal{L}_3$ in the case $t=3$ in the above succession rule
$\Omega_t$.

\begin{description}
    \item \begin{description}
                \item $\mathcal{L}_1=<3>;$
          \end{description}
    \item $L_2^1=3\circ s(3,3)=3 \circ <3,5,4>=<33,35,34>$, then
    \item \begin{description}
                \item $\mathcal{L}_2=<33,35,34>;$
          \end{description}
    \item $L_3^1=33\circ s(3,3)=33\circ <3,5,4>=<333,335,334>$;
    \item $L_3^2=35\circ s(5,4)=35\circ <4,5,6,7,3>=<354,355,356,357,353>$;
    \item $L_3^3=34\circ s(4,3)=34\circ <3,6,5,4>=<343,346,345,344>$, then
    \item \begin{description}
                \item
                $$\mathcal{L}_3=<333,335,334,354,355,356,357,353,343,346,345,344>.$$
          \end{description}
\end{description}
If $t=2$, then we find the succession rule $\Omega_C$ for Catalan
numbers, enumerating, among other things, the binary trees. In
\cite{V1} the author proposes a constant time algorithm for
generating binary trees Gray codes. We note that our procedure,
combined with the results of Section \ref{Dyck}, is an alternative
approach for this aim.
\begin{flushright}$\square$\end{flushright}\phantom\\

\vspace{.3cm}\noindent \textbf{Example 3.} Another particular case
of $\Omega_A$ is the following family:
$$
\Omega_r: \left\{
\begin{array}{l}
   (r)\\
   (k)\rightsquigarrow(r)(r+1)\ldots(k)(k+1)^{r-1}\ ,\\
\end{array}
\right.
$$
with $r\geq 2$. They satisfy the stability property, too, with
$i=1$, $j=2$, $c_1=r$ and $c_2=r+1$. If $r=3$, then $\Omega_r$ is
the well-known succession rule defining the sequence of Schr\"oder
numbers. The following rule $\Omega_s$ also codes the construction
of Schr\"oder paths, 2-colored parallelogram polyominoes, (4231,
4132)-pattern avoiding permutations, (3142, 2413)-pattern avoiding
permutations \cite{BDPP2,We1,We2} (these latter patterns are also
considered in \cite{BBL} for pattern matching decision problem for
permutations).
$$
\Omega_s: \left\{
\begin{array}{l}
   (2)\\
   (k)\rightsquigarrow(3)(4)\ldots(k)(k+1)^2\\
\end{array}
\right.
$$
In this case it is $c_1=3$, $c_2=4$ and the associated shifted
production is:
$$
\left\{
\begin{array}{l}
   s(k,3)=<3,(k+1)_2,(k+1)_1,k,k-1,\ldots,5,4>\\
   \\
   s(k,4)=<4,5,\ldots,k,(k+1)_1,(k+1)_2,3>\ ,\\
\end{array}
\right.
$$
where the indexes differentiate labels with the same value. Note
that $s(k_2,*)=s(k_1,*)$ ($*=3\ \mbox{or}\ 4$). The construction
of the list $\mathcal{L}_3$ is:

\begin{description}
    \item \begin{description}
                \item $\mathcal{L}_1=<3>;$
          \end{description}
    \item $L_2^1=3\circ s(3,3)=<33,34_2,34_1>$, then
    \item \begin{description}
                \item $\mathcal{L}_2=<33,34_2,34_1>;$
          \end{description}
    \item $L_3^1=33\circ s(3,3)=<333,334_2,334_1>$;
    \item $L_3^2=34_2\circ s(4_2,4_1)=<34_24_1,34_25_1,34_25_2,34_23>$;
    \item $L_3^3=34\circ s(4,3)=<343,345_2,345_1,344>$, then
    \item \begin{description}
                \item
                $$\mathcal{L}_3=<333,334_2,334_1,34_24_1,34_25_1,34_25_2,
                34_23,343,345_2,345_1,344>.$$
          \end{description}
\end{description}
\begin{flushright}$\square$\end{flushright}\phantom\\

\vspace{.3cm}\noindent \textbf{Example 4.} Succession rules of
kind:

$$
\Omega_B: \left\{
\begin{array}{l}
   (r)\\
   (k)\rightsquigarrow(b)^l(a)(a+1)\ldots(k)(k+d_1)\ldots(k+d_m)\\
   (k)\rightsquigarrow(b)^k \qquad \qquad\qquad \qquad\qquad \qquad if \ (k<a)\wedge(k\leq l)\\
   (k)\rightsquigarrow(b)^{(k-1)}(a) \qquad \qquad\qquad \qquad\qquad if \ k<a \ \ , \\
\end{array}
\right.
$$
with $l\geq 2$, $b<a$, $m=a-l-1$, satisfy the stability property
with $i=1$, $j=2$ and, denoting $b^l=b_1b_2\ldots b_l$, $c_1=b_1$,
$c_2=b_2$. A well-known particular case is
$$
\Omega_{GD}: \left\{
\begin{array}{l}
   (2)\\
   (2)\rightsquigarrow(3)(3)\\
   (3)\rightsquigarrow(3)(3)(4)\\
   (k)\rightsquigarrow(3)^2(4)\ldots(k)(k+1)\\
\end{array}
\right.
$$
which encodes a construction for Gran Dyck paths \cite{PPR}. The
generalized shifted production associated is
$$
\left\{
\begin{array}{l}
   s(k,3_1)=<3_1,k+1,k,\ldots,4,3_2>\\
   \\
   s(k,3_2)=<3_2,4,\ldots,k,k+1,3_1>\ .\\
\end{array}
\right.
$$
The list $\mathcal{L}_3$ is obtained as follows:

\begin{description}
    \item \begin{description}
                \item $\mathcal{L}_1=<2>;$
          \end{description}
    \item $L_2^1=2\circ s(2,3_1)=<23_1,23_2>$, then
    \item \begin{description}
                \item $\mathcal{L}_2=<23_1,23_2>;$
          \end{description}
    \item $L_3^1=23_1\circ s(3_1,3_1)=<23_13_1,23_14,23_13_2>$;
    \item $L_3^2=23_2\circ s(3_2,3_2)=<23_23_2,23_24,23_23_1>$, then
    \item \begin{description}
                \item
                $$\mathcal{L}_3=<23_13_1,23_14,23_13_2,23_23_2,23_24,23_23_1>.$$
          \end{description}
\end{description}
\begin{flushright}$\square$\end{flushright}\phantom\\

\vspace{.3cm}\noindent \textbf{Example 5.} It is possible to find
some examples among the transcendental succession rules which are
stable or not. The classical rule defining the factorial numbers,
which describes the construction of  the permutations of length
$n$ by inserting the element $n$ in any active site of any
permutation of length $n-1$, is not stable (its production is:
$(k)\rightsquigarrow (k+1)^k$). On the contrary, the following one
$\Omega_p$, defining the same sequence, is stable:
$$
\Omega_p= \left\{
\begin{array}{l}
(2)\\
(2k)\rightsquigarrow(2)(4)(6)\ldots(2k)(2k+2)^k\ .
\end{array}
\right.
$$
Stability property is satisfied since each label $(2k)$ generates
in the first two positions labels (2) and (4). The associated
generalized shifted production is:
$$
\left\{
\begin{array}{l}
   s(2k,2)=<2,(2k+2)_k,(2k+2)_{k-1},\ldots,(2k+2)_1,2k,2k-2,\ldots,4>\\
   \\
   s(2k,4)=<4,6,\ldots,2k-2,2k,(2k+2)_1,(2k+2)_2,\ldots,(2k+2)_k,2>\ \ ,\\
\end{array}
\right.
$$
where the indexes are useful to distinguish different labels but
with the same value.
 In order to illustrate the
combinatorial placement of $\Omega_p$ we propose a probably new
ECO construction for the permutations which can be described by
this rule. Let $\pi= \pi_1\pi_2\ldots \pi_n$ be a permutation of
$S_n$, we define an operator $\vartheta: S_n \longrightarrow
2^{S_{n+1}}$ (the power set of $S_{n+1}$) working as follows
($n\geq 1$):
\begin{itemize}
    \item let $\pi_1=k$, then $\vartheta$ generates $2k$ permutations $\pi'\in
    S_{n+1}$ which are indicated by $\pi'^{(i)}$, with $i=1,2,\ldots,2k$;
    \item the entries of $\pi'^{(i)}$ are:
    \begin{enumerate}
        \item if $i=1,2,\ldots,k$, then:
        \begin{description}
            \item[$\circ$] $\pi'^{(i)}_1=i$;
            \item[$\circ$] the other entries are the same of $\pi$
            where the entry $i$ is replaced by $n+1$.
        \end{description}
        \item if $i=k+1,k+2,\ldots,2k$, then:
        \begin{description}
            \item[$\circ$] $\pi'^{(i)}_1\pi'^{(i)}_2=(\pi_1+1)j$, where $j=1,2,\ldots,k$;
            \item[$\circ$] the other entries are obtained as follows:
            \begin{itemize}
                \item If $\pi_1\neq n$, then let $\rho$ be the
                sequence, with length $n-1$, obtained by $\pi$
                deleting $\pi_1$ after it has been interchanged
                with $\pi_1+1$. The remaining entries of $\pi'^{(i)}$
                are the same of $\rho$ where the entry $j$ is
                replaced by $n+1$.
                \item If $\pi_1=n$, then let $\rho$ be the sequence
                obtained from $\pi$ by deleting $\pi_1$. The
                remaining entries of $\pi'^{(i)}$ are the same of
                $\rho$ where the entry $j$ is replaced by $n$.
            \end{itemize}
        \end{description}
    \end{enumerate}
\end{itemize}

\noindent \textbf{Remark:} permutations $\pi'^{(i)}$ with
$i=1,2,\ldots,k$ start with an ascent, while permutations
$\pi'^{(i)}$ with $i=k+1,k+2,\ldots,2k$ start with a descent.
\\ \par
It can be easily proved that if $\pi'\in S_{n+1}$, then there
exists a unique $\pi\in S_n$ such that $\pi'\in \vartheta(\pi)$
($n\geq 1$), then operator $\vartheta$ satisfies Proposition 2.1
of \cite{BDPP1}, which ensures that the family of sets
$\{\vartheta(\pi):\pi\in S_n\}$ is a partition of $S_{n+1}$, so
that $\vartheta$ provides a recursive construction of the
permutations of $S=\bigcup S_n$.

\begin{figure}
\begin{center}
\includegraphics[scale=.5]{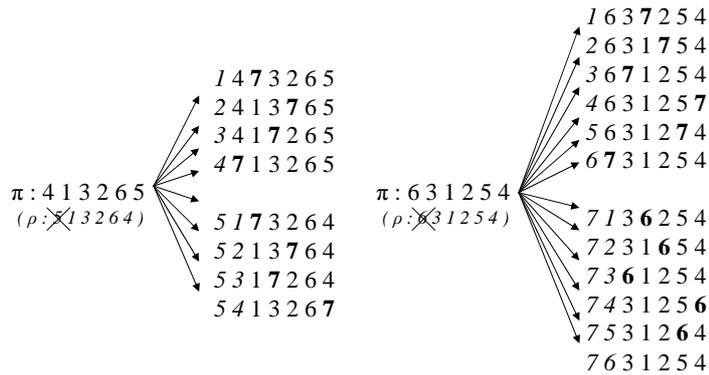}
\end{center}
\caption{The action of $\vartheta$ on two different permutations
of $S_6$.}\label{fig:permutazioni}
\end{figure}

In Figure \ref{fig:permutazioni} the action of $\vartheta$ on two
different permutations of $S_6$ (the first one starting with an
entry different from $n=6$) is illustrated. Permutations
$\pi'^{(i)}$, $i=1,2,\ldots,2k$ generated by $\pi$ by means
$\vartheta$ are listed from the top to the bottom, being
$\pi'^{(1)}$ at the top.

\section{Conclusions and further developments}

It is possible to find a lot of succession rules satisfying the
stability property, but we are interested to the rules having some
combinatorial relevance, as the ones presented in the above
examples. In this way, with our procedure we are able to give a
Gray code for the words (i.e. the paths whose nodes are the labels
in the generating tree) encoding combinatorial Gray structures,
i.e. those structures whose exhaustive generation can be described
by a rule satisfying the stability property, which is not, as we
have seen, an infrequent property.

Clearly, it would be better to have a Gray code for the objects
instead of their encodes. Nevertheless, as we stated in Section
\ref{intro}, our procedure generates a Gray code which is not
related to the nature of a particular class of combinatorial
objects. Moreover, in some case it could be possible to translate
the word of labels (the path in the generating tree) into the
corresponding object. A further effort in this sense could be the
research of algorithms for this translation in order to generalize
the approach of Section \ref{Dyck} for Dyck paths. For this aim
the ECO method can be useful, since by means of it each code is
associated to a single object of the structure.

From the above examples it is possible to argue that the stability
property of a succession rule does not depend on its "structural
properties`` , which have been discussed by the authors in
\cite{BETAL}. In the light of this fact, it is reasonable to ask
if a stable succession rule can be considered as the
representative, say \emph{standard form}, of a set of rules which
are all equivalent to it (two rules are said equivalent if they
define the same number sequence \cite{BDPR}). This is to say that
the equivalence problem for succession rules could be amplified
with respect to the investigation conducted in \cite{BDPR} where
the authors analyze the equivalence problem for some different
kinds of rules: is it suitable the research of the set of rules
equivalent to a stable succession rule?

Moreover, it is evident that it is not the sequence defined by the
rule that induces it to be stable or not: factorial number
sequence can be defined by a stable or not stable rule, as showed
in Example 5. Consequently, a problem which naturally arises from
this note is the existence of a succession rule with the stability
property for any given number sequence. A first concerning
question could be the following (to the authors knowledge the
answer is open): is there a stable rule defining Motzkin numbers?



\end{document}